# DISCUSSION: ONE-STEP SPARSE ESTIMATES IN NONCONCAVE PENALIZED LIKELIHOOD MODELS: WHO CARES IF IT IS A WHITE CAT OR A BLACK CAT?[1]

By Xiao-Li Meng

*Harvard University*

**1. An insider's minor comments.** Section 2.3 seems to be the reason that I am a discussant. There, it was first stated that the proposed LLA algorithm is an instance of the MM algorithm, as termed by Lange, Hunter and Yang [8]. Then it was shown, under certain conditions, that it is also an EM algorithm. My initial reaction was "hmmm, the authors' reading of Lange, Hunter and Yang [8] must have ceased before reaching its discussions," because a more general "MM = EM" result using the same Laplace transform technique was the base for a key inquiry of Meng [10], a discussion of Lange, Hunter and Yang [8]. Upon a more careful reading, I realized that the authors' construction, though mathematically equivalent to mine, gives a different interpretation to the constructed missing data/latent variable. This is rather interesting, especially if my initial reaction was correct.

For the current paper, this "MM = EM" result appears to be of minor interest, especially because its potential benefit is not explored in the paper. Instead, the only punch line seems to be a logically unsubstantiated one: "Thus, Theorem 3 also indicates that MM algorithms are more flexible than EM algorithms." To clarify these issues, which seems to be what I have been asked to do, I'll focus my discussion first on this algorithmic connection and qualify the statements that I just made. I will then pose some questions going beyond those algorithmic and computational considerations. (The two parts are connected via the "cat" title, for those who are patient enough to read both.)

To do so most effectively, let me invoke an author's privilege to reproduce Section 2 of Meng [10] in its entirety, to map out its mathematical connections with Section 2.3 of the current paper. The acronym "SM" below was

---

Received November 2007; revised November 2007.

[1]Supported in part by several NSF grants.

my suggestion because the general recipe Lange, Hunter and Yang [8] put forward consists of a "Surrogate step" and a "Maximization/Minimization step"; in their rejoinder, Lange, Hunter and Yang [8] coined the term "MM," partially in fear of the association of "SM" with "S&M." The "alphabet soup" phrase was used by Lange, Hunter and Yang [8] to describe the collection of acronyms in the EM literature; indeed, for readers who enjoy collecting acronyms, the "GAECM" below is another treasure to hunt!

**2. Is SM just EM?** Of course, the new alphabet soup will not be a really new delight if it is just the old soup presented in a new, perhaps larger, bowl. Could it be that SM, though apparently more general, is just a disguised or beautified version of EM? The answer is not completely obvious, especially if one starts the comparison with the most obvious construction of the surrogate function via linear minorization/majorization. As in the authors' equation (3.1) (page 9), assume our log-likelihood function $L(\theta|y)$ can be written as

$$(2.1) \qquad L(\theta|y) = f_y(\theta) - g_y(\theta), \qquad \theta \in \Theta \subset R^1,$$

where both $f_y$ and $g_y$ are concave functions and without loss of generality [when $|g_y(0)| < \infty$] we assume $g_y(0) = 0$ for all $y$. Now suppose $e^{-g_y(\theta)}$ is the moment generating function of a conditional density $h(z|y)$, namely,

$$(2.2) \qquad e^{-g_y(\theta)} = \int e^{\theta z} h(z|y) \mu(dz), \qquad \theta \in \Theta.$$

Then if we augment $p(y|\theta) = e^{L(\theta|y)}$ to

$$(2.3) \qquad p(z|y,\theta)p(y|\theta) \equiv [e^{\theta z + g_y(\theta)} h(z|y)][e^{L(\theta|y)}] = e^{f_y(\theta) + \theta z} h(z|y),$$

we have, for the standard EM construction,

$$(2.4) \qquad Q(\theta|\theta^{(t)}) = f_y(\theta) + \theta \mathrm{E}(Z|y,\theta^{(t)}) + \mathrm{E}[\log h(Z|y)|y,\theta^{(t)}].$$

But this is equivalent to the proposed linear minorization surrogate function

$$(2.5) \qquad Q(\theta|\theta^{(t)}) = f_y(\theta) - g'_y(\theta^{(t)})(\theta - \theta^{(t)}),$$

because $\mathrm{E}(Z|y,\theta) = -g'_y(\theta)$ from differentiating both sides of (2.2). Incidentally, by differentiating both sides of (2.2) twice, we have $g''_y(\theta) = -\mathrm{V}(Z|y,\theta) \leq 0$, and thus the concavity of $g(\theta)$ is a necessary condition for this EM construction to be possible. [For multivariate $\theta$ we can construct the missing data $Z$ with the same dimension and replace $\theta z$ in (2.2) with $\theta^{\top} z$.]

Although *this* EM construction is not always possible (e.g., when $e^{-g_y(\theta)}$ may not be a moment generating function), and even when it is possible it requires more brain power than the linear minorization method, it nevertheless suggests that a large class of SM algorithms based on (2.5) are also EM algorithms with augmentation $p(z,y|\theta)$ of (2.3).

So the question is, given $Q(\theta|\phi)$ from a particular SM construction, how do we know if there is a corresponding EM construction, regardless of how convoluted the latter might be? The practical relevance of this theoretically interesting question is that, if the EM class is as rich as the SM class, then the value of the new SM formulation is in providing a set of new tools for creative EM-type implementation. However, if the SM class is richer than the EM class, then it provides hope for solving problems that are difficult or even impossible to solve within the entire GAECM framework.

To apply the above result to the setting of authors' (2.1), we simply set $\theta = \beta$ and let $f_y(\theta) = \sum_{i=1}^{n} \ell_i(\beta)$, $g_y(\theta) = n \sum_{j=1}^{p} p_{\lambda_j}(\beta_j)$. Then the authors' (2.1) can be expressed as my (2.1) if we replace the $\theta \equiv \beta$ in $g_y(\theta)$ by $|\beta| \equiv (|\beta_1|, \ldots, |\beta_p|)^{\top}$. The authors' equivalence result then follows if we make the same replacement of $\theta$ by $|\theta|$ in my augmentation (2.3) above, that is, if we augment $p(y|\theta) = e^{L(\theta|y)}$ to

$$(1.1) \quad p(z|y,\theta) p(y|\theta) \equiv [e^{|\theta|^{\top} z + g_y(|\theta|)} h(z|y)][e^{L(\theta|y)}] = e^{f_y(\theta) + |\theta|^{\top} z} h(z|y).$$

A key difference between the authors' construction and mine is that in the authors' setting, the *augmented variable*—or *auxiliary variable*, as known in statistical physics (but AV in either case!)—$\tau$ is a hyper-parameter, while in my construction, the AV $z$ represents missing data. Although I believe that the authors' reciprocal transformation, that is, taking $\tau \propto z^{-1}$, is unnecessary and only mystifies the construction, the authors' hyper-parameter interpretation further helps to emphasize the key question raised in the last paragraph of the quoted Section 2 above. That is, given an MM (a.k.a. SM) algorithm, how do we know it is *not* an EM algorithm, considering there are infinitely many ways to construct AVs? Just because one particular construction fails to reproduce an MM algorithm as an EM algorithm, as in authors' Theorem 3, logically it says nothing about whether this MM algorithm can be reproduced as an EM algorithm by a different AV construction.

In fact, the identification of whether an MM algorithm is also an EM algorithm (the reverse of course is always true) is an unexpectedly difficult problem. In Meng [10], I gave a sufficient and necessary condition for an MM to be an EM. Verifying this condition turns out to be the same as determining when a *normalized yoke* is also an *expected likelihood yoke*, a very challenging problem in differential geometry. Any reader who is curious about the notion of *yoke*, or wants to take on this challenge, is invited to spend an intensive evening with me via Meng [10], with, of course, guests such as Barndorff–Nielsen [2], Barndorff–Nielsen and Jupp [3] and Blæsild [4]—there could be an AOS-worthy paper here!

I surmise most readers of the current paper would prefer to spend an evening with someone else. Indeed, I am a bit curious myself about the connection of this equivalence result to the rest of the paper. "Who cares if it is a white cat or black cat, as long it catches mice," so goes a famous quotation attributed to Deng Xiaoping (allegedly the statement got him into trouble before he regained his power in China). So who cares if it is an MM or EM, as long as it finds the right solution? One answer is that once we know an algorithm is an EM algorithm, we can take advantage of the vast literature on the EM algorithm to study its theoretical properties or improve its speed by applying techniques designed for EM-type of algorithms. I can certainly think of the benefit of deriving its theoretical rate of convergence along the lines of Meng and Rubin [11, 12] or Meng [9], or trying to speed it up by attempting the recipes given in Meng and van Dyk [13, 14], either of which could lead to useful insights or improvement of LLA. Seeing no such discussions in this paper, I was almost tempted to self-indulge the possibility that this MM-EM equivalence issue was included because the authors wanted me to be a discussant (circumstantial evidences: [10] is the only general investigation of the equivalence result, to the best of my knowledge). Of course, my ego-checking side has to wonder: where is the beef (or mouse!)?

**2. An outsider's major questions.** Having done my insider's job but failed to identify a mouse, let me now try to catch a bird by going outside! As an outsider (to penalized likelihood), I gather I can enjoy the freedom of asking naive questions and rely on the authors' expertise to correct any "mis-information" my inquiries might have created for the readers.

**Q.1.** *Should we care about the difference between penalized likelihood methods and Bayesian methods?*

To many Bayesians, and perhaps some semi-Bayesian or even non-Bayesians, using penalized likelihood is an attempt of enjoying the Bayesian fruits without paying the B-club fee. Although some pena-likelihoodists might be offended by this analogy, many of us perhaps are of the "who cares?" opinion here again: who cares if it has a PL-label or B-label, as long as it gets the job done? Indeed, I had been one of those until I gave the question a bit more thought when I was studying the literature on wavelet thresholding, a form of variable selection.

The central question here is: are the two methods as general and flexible as each other, so the difference is only in perception not in operation? On the surface they appear to be: just equate the penalty function term with the log of the prior density—the fact that Bayesian methods allow improper prior densities seems to make this equivalence almost completely general. However, in the Bayesian framework, once a posterior distribution is obtained, there is more than one way to utilize it; in the pena-likelihood world, maximizing the likelihood seems to be the only general recipe.

A main focus of the current paper is about finding efficient and practical algorithms for maximizing the penalized log-likelihood in the form of authors' (2.1). Various approximations were made to the penalty function to overcome technical difficulties. When confined in the pena-likelihood framework, this seems to be the only game in town, as least from an outsider's perspective. When recast in the Bayesian framework, at least three additional considerations come to my mind, however.

First, considering the penalty function as the log of the prior density, it seems to me more fruitful to regard each approximation in its own right as a prior specification, rather than an approximation to an "ideal" penalty function or log-prior. Evidently, LQA can be viewed as using a locally-adaptive

Gaussian prior, and LLA a locally-adaptive Laplace prior. Viewed in this way, it is almost inevitable to wonder what will happen if we adopt a locally-adaptive $t$-distribution as our prior or some other families of distributions? Even for people who simply do not want the B-club membership, these links can potentially provide new approximations to penalty functions that are otherwise hard to come by. For example, using mixtures as prior distributions is a common practice in the Bayesian world, but it may not come as naturally when one is in the penalty function mood, though the authors' (2.2) hints at that direction.

Second, a key advantage of using an $L^1$ penalty appears to be its ability to ensure sparse representation. From a Bayesian perspective, a sparse representation means that a priori we should put nontrivial mass at $\beta_j = 0$ for any $j$ (since a priori we do not know which one is zero). This leads naturally to the mixture prior specifications, such as those used in Abramovich, Sapatinas and Silverman [1] and Johnston and Silverman [5]. The need for nontrivial mass at zero seems to offer an alternative insight for the nondifferentiability at zero of the SCAD penalty family, as the authors discussed. Indeed, only thinking in terms of prior specifications allowed me, as an outsider, to understand intuitively why the optimal penalty function has to be nondifferentiable at the origin.

Third, once a posterior distribution is in place, we can consider several point-estimate summaries other than its mode (i.e., corresponding to MLE). Indeed, as demonstrated in Abramovich, Sapatinas and Silverman [1], if we use the posterior *median*, then we can also ensure sparse representation under the mixture priors mentioned above. To see this most clearly in general, suppose our prior for a (scalar) regression coefficient $\beta_j$ is a mixture of an atom at $\beta_j = 0$ and a continuous component, whose cumulative distribution function (CDF) is given by $F_j(\beta_j)$, $j = 1, \ldots, p$. Under the assumption of a priori independence, the joint prior CDF for $\{\beta_1, \ldots, \beta_p\}$ is given by

$$(2.1) \qquad F(\beta_1, \ldots, \beta_p) = \prod_{i=1}^{p} [\pi_j 1_{(\beta_j \geq 0)} + (1 - \pi_j) F_j(\beta_j)],$$

where $\pi_j$ is the prior mass on $\beta_j = 0$, governing our desired degree of sparsity (e.g., setting $\pi_j = 0.9$ for all $j$ means that we a priori believe that only about 10% of the $\beta$ coefficients should be significant).

Under this setting, it is easy to verify that the *marginal* posterior distribution for any $\beta_j$ is also a mixture of an atom at zero and a continuous component:

$$(2.2) \qquad F_j(\beta_j | \mathbf{y}) = \pi_{\mathbf{y},j} 1_{(\beta_j \geq 0)} + (1 - \pi_{\mathbf{y},j}) F_{c,j}(\beta_j | \mathbf{y}).$$

The exact expressions for $\pi_{\mathbf{y},j}$ and $F_{c,j}(\beta_j | \mathbf{y})$ follow directly from the Bayes theorem, and they are the simplest when the likelihood $L(\beta_1, \ldots, \beta_p)$ factors

into $\prod_j L_j(\beta_j)$, as in the wavelet applications of Abramovich, Sapatinas and Silverman [1] and Johnston and Silverman [5], in which case the $\{\beta_1, \ldots, \beta_p\}$ are a posteriori independent. But regardless of this independence, as long as we decide to threshold each coefficient separately (in contrast to "block thresholding"), which seems to be the case for the current paper, we only need to deal with the marginal distribution as given in (2.2). We therefore, for notation simplicity, will drop the subscript $j$ in the following discussion. (For a "true" Bayesian such a marginal approach can also be acceptable when the dependence is not too strong and "joint thresholding" is computationally too expensive to be practical. But even in such cases the Bayesian thinking clearly helps us to understand what corners we have cut, thereby constructively pointing us to directions for improvement when more efficient computational tools become available.)

Given the setting above, it is trivial to see that if the posterior mass $\pi_\mathbf{y} \equiv \Pr(\beta = 0 | \mathbf{y}) \geq 1/2$, then $\Pr(\beta \leq 0 | \mathbf{y}) \geq 1/2$ and $\Pr(\beta < 0) \leq 1/2$. Recall that the median of a (right-continuous) CDF $F(x)$ is the (or any) value $M$ such that

$$(2.3) \qquad F(M) \geq \tfrac{1}{2} \quad \text{but} \quad F(M_-) \equiv \lim_{x \uparrow M} F(x) \leq \tfrac{1}{2}.$$

Consequently, whenever $\pi_\mathbf{y} \geq 1/2$, the posterior median of $\beta$, denoted by $\mathrm{Med}(\beta | \mathbf{y})$, must be exactly zero, provided that the continuous component $F_c$ has no "flat" region, a condition we assume. [Without this assumption, $\mathrm{Med}(\beta | \mathbf{y})$ may not be unique, a technical complication that is of little interest for almost any practical purposes.] In other words, the median procedure will exclude a coefficient $\beta_j$ if the posterior lends stronger (or no weaker) support to $\beta_j = 0$ than to $\beta_j \neq 0$, a very intuitive and appealing requirement for sparsity to occur in our inference.

Given that $\pi_\mathbf{y} \geq 1/2$ is a sufficient condition for declaring sparsity under the posterior median approach, it is natural to ask if it is also a necessary condition. That is, is it possible for $\mathrm{Med}(\beta | \mathbf{y})$ to be exactly zero when $\pi_\mathbf{y} < 1/2$? The answer is positive, as demonstrated in Abramovich, Sapatinas and Silverman [1]. This might appear to be somewhat counterintuitive at first sight, because $\pi_\mathbf{y} < 1/2$ seems to imply that the hypothesis $\beta = 0$ is less supported by the data than the hypothesis $\beta \neq 0$ is. There is a subtle distinction here, however. A lack of support for the atom component of the mixture distribution is not the same as a lack of support for the value of $\beta$ to be zero. The continuous part of the posterior distribution, namely, $F_c(\beta | \mathbf{y})$, can also provide evidence *for* $\beta = 0$, even if it was "designed" for inference for $\beta$ away from zero. This perhaps can be best seen in the extreme case where we set $\pi = 0$, that is, we do not believe sparsity a priori and, hence, we do not even allow the atom component to show up in our model. Nevertheless, our posterior inference can still provide good evidence for $\beta = 0$, or more precisely, for $\beta$ being in a small neighborhood of $\beta = 0$, if, for example, our posterior density ends up like a normal density with mean (near) zero.

This observation suggests that even when $\pi_{\mathbf{y}} < 1/2$, there should be a threshold ($< 1/2$) for $\pi_{\mathbf{y}}$ above which $\mathrm{Med}(\beta|\mathbf{y})$ would still be zero because the additional evidence for $\beta = 0$ from the continuous component $F_c(\beta|\mathbf{y})$ is able to compensate the degree of failing of reaching $\pi_{\mathbf{y}} \geq 1/2$ by the atom component. To see this clearly, let $\mathrm{Med}_c(\beta|\mathbf{y})$ be the posterior median from $F_c(\beta|\mathbf{y})$. If $\mathrm{Med}_c(\beta|\mathbf{y}) > 0$, then for any $\beta < 0$,

$$(2.4) \qquad F(\beta|\mathbf{y}) = (1 - \pi_{\mathbf{y}})F_c(\beta|\mathbf{y}) < \frac{1 - \pi_{\mathbf{y}}}{2} \leq \frac{1}{2}.$$

Hence, $\mathrm{Med}(\beta|\mathbf{y})$ must be nonnegative. This could let us to declare that $\mathrm{Med}(\beta|\mathbf{y})$ must be the solution of

$$(2.5) \qquad F(\beta|\mathbf{y}) = \pi_{\mathbf{y}} + (1 - \pi_{\mathbf{y}})F_c(\beta|\mathbf{y}) = \tfrac{1}{2},$$

solving which would yield $\mathrm{Med}(\beta|\mathbf{y}) = F_c^{-1}(\frac{1 - O_{\mathbf{y}}}{2}|\mathbf{y})$, where $O_{\mathbf{y}} = \pi_{\mathbf{y}}/(1 - \pi_{\mathbf{y}})$ is the posterior odds for $\beta = 0$. One has to realize, however, that although equation (2.5) always has a unique solution whenever $O_{\mathbf{y}} < 1$, this solution is nonnegative if and only if $O_{\mathbf{y}} \leq 1 - 2F_c(0|\mathbf{y})$, which is a nonempty condition because $F_c(0|\mathbf{y}) \leq 1/2$ under our assumption that $\mathrm{Med}_c(\beta|\mathbf{y}) > 0$. When this condition fails, clearly $F(0|\mathbf{y}) > 1/2$, and hence, $\mathrm{Med}(\beta|\mathbf{y}) = 0$, because $F(0_-|\mathbf{y}) \leq 1/2$, a consequence of (2.4).

Similarly, when $\mathrm{Med}_c(\beta|\mathbf{y}) < 0$, then for any $\beta > 0$,

$$(2.6) \qquad F(\beta_-|\mathbf{y}) = \pi_{\mathbf{y}} + (1 - \pi_{\mathbf{y}})F_c(\beta|\mathbf{y}) > \pi_{\mathbf{y}} + \tfrac{1}{2}(1 - \pi_{\mathbf{y}}) \geq \tfrac{1}{2}.$$

Therefore, $\mathrm{Med}(\beta|\mathbf{y})$ must be nonpositive. This means that $\mathrm{Med}(\beta|\mathbf{y})$ is the solution of $(1 - \pi_{\mathbf{y}})F_c(\beta|\mathbf{y}) = 1/2$, *when* the solution $F_c^{-1}(\frac{1 + O_{\mathbf{y}}}{2}|\mathbf{y}) \leq 0$, which holds if and only if $O_{\mathbf{y}} \leq 2F_c(0|\mathbf{y}) - 1$. When this condition fails, $F(0_-|\mathbf{y}) < 1/2$, and hence, $\mathrm{Med}(\beta|\mathbf{y}) = 0$ because (2.6) implies $F(0) \geq 1/2$.

Finally, when $\mathrm{Med}_c(\beta|\mathbf{y}) = 0$,

$$F(0_-|\mathbf{y}) = \frac{1 - \pi_{\mathbf{y}}}{2} \leq \frac{1}{2} \quad \text{and} \quad F(0|\mathbf{y}) = \frac{1 + \pi_{\mathbf{y}}}{2} \geq \frac{1}{2},$$

hence, $\mathrm{Med}(\beta|\mathbf{y}) = 0$. Summarizing these derivations yields the following general result, which provides a theoretically more appealing form to gain insights than detailed expressions derived for specific models, such as the formulae in Abramovich, Sapatinas and Silverman [1] under normality.

LEMMA 1. *Assume the continuous component $F_c(\beta|\mathbf{y})$ in the posterior mixture is strictly monotone, and let $S_{c,\mathbf{y}} = \mathrm{Sign}\{\mathrm{Med}_c(\beta|\mathbf{y})\}$, the sign of $\mathrm{Med}_c(\beta|\mathbf{y})$. Then the posterior median of $\beta$ has the following analytic expression, where $O_{\mathbf{y}}$ is the posterior odds for $\beta = 0$:*

$$(2.7) \quad \mathrm{Med}(\beta|\mathbf{y}) = \begin{cases} F_c^{-1}\left(\dfrac{1 - S_{c,\mathbf{y}}O_{\mathbf{y}}}{2}\Big|\mathbf{y}\right), & \text{if } O_{\mathbf{y}} < |1 - 2F_c(0|\mathbf{y})|, \\ 0, & \text{if } O_{\mathbf{y}} \geq |1 - 2F_c(0|\mathbf{y})|. \end{cases}$$

This result, of course, is elementary (so no credit is claimed!), but the insights deriving from it are less so. First, the quantity $\Delta_c(\mathbf{y}) = |1 - 2F_c(0|\mathbf{y})|$ provides the threshold we were looking for, that is, the smallest value of $O_{\mathbf{y}}$ such that if $O_{\mathbf{y}}$ takes any value below it, then the evidence for $\beta = 0$ from the continuous component can no longer compensate for the fact that $O_{\mathbf{y}} < 1$ [note that $O_{\mathbf{y}} \geq \Delta_c(\mathbf{y})$ holds trivially when $O_{\mathbf{y}} \geq 1$]. Because

$$\Delta_c(\mathbf{y}) = |\Pr_c(\beta > 0|\mathbf{y}) - \Pr_c(\beta \leq 0|\mathbf{y})|,$$

where the probability calculation $\Pr_c(\cdot|\mathbf{y})$ is with respect to $F_c(\cdot|\mathbf{y})$, we see clearly that the closer $\mathrm{Med}_c(\beta|\mathbf{y})$ is to zero—and hence the more "borrowed" evidence for supporting $\beta = 0$—the more tolerant the median approach becomes in setting the lower bound on the permissible value of $O_{\mathbf{y}}$ for hard-thresholding $\beta$. Indeed, in the extreme case when $\mathrm{Med}_c(\beta|\mathbf{y})$ is exactly zero, the posterior median approach will set $\beta = 0$ regardless of the value of $O_{\mathbf{y}}$. This makes sense because $\mathrm{Med}_c(\beta|\mathbf{y}) = 0$ is the strongest evidence for $\beta = 0$ from the continuous component under the $L^1$ norm (and when only point estimators are considered).

Second, even when $O_{\mathbf{y}} < \Delta_c(\mathbf{y})$, there is still some evidence for $\beta = 0$, and the posterior median approach directly takes into account the strength of this evidence by shrinking $\mathrm{Med}_c(\beta|\mathbf{y}) = F_c^{-1}(\frac{1}{2}|\mathbf{y})$ toward zero according to the value of $O_{\mathbf{y}}$. By Lemma 1, $\mathrm{Med}(\beta|\mathbf{y})$ is always closer to zero than $\mathrm{Med}_c(\beta|\mathbf{y})$, and how much closer is directly determined by $O_{\mathbf{y}}$. Third, this shrinkage is actually a double shrinkage, representing a sensible form of "soft thresholding." This is perhaps most clearly seen when $F_c(\beta|\mathbf{y})$ is given by a symmetric location-scale family (e.g., normal): $G(\frac{\beta - \mathrm{Med}_c(\beta|\mathbf{y})}{\sigma_{c,\mathbf{y}}})$, where $\sigma_{c,\mathbf{y}}$ is the posterior scale, measuring the uncertainty in the continuous component of the posterior distribution [and in $\mathrm{Med}_c(\beta|\mathbf{y})$] for inferring $\beta$. In such cases, Lemma 1 implies that, when $O_{\mathbf{y}} < \Delta_c(\mathbf{y})$ [and hence $\mathrm{Med}_c(\beta|\mathbf{y}) \neq 0$],

$$\mathrm{Med}(\beta|\mathbf{y}) = \begin{cases} \mathrm{Med}_c(\beta|\mathbf{y}) - \sigma_{c,\mathbf{y}} G^{-1}\left(\dfrac{1 + O_{\mathbf{y}}}{2}\right), & \text{if } \mathrm{Med}_c(\beta|\mathbf{y}) > 0; \\[2mm] \mathrm{Med}_c(\beta|\mathbf{y}) + \sigma_{c,\mathbf{y}} G^{-1}\left(\dfrac{1 + O_{\mathbf{y}}}{2}\right), & \text{if } \mathrm{Med}_c(\beta|\mathbf{y}) < 0. \end{cases}$$

(2.8)

The "double shrinkage" is seen here because $\mathrm{Med}_c(\beta|\mathbf{y})$ is often a shrinkage estimator by itself (e.g., in the normal case and when the prior mean for $F_c$ is set to zero, as typical in applications). But as long as $O_{\mathbf{y}} > 0$ and $\sigma_{c,\mathbf{y}} > 0$ (always true in practice), $\mathrm{Med}(\beta|\mathbf{y})$ further shrinks toward zero. And the larger $\sigma_{c,\mathbf{y}}$ or $O_{\mathbf{y}}$, the greater the shrinkage. This makes perfect sense, as larger $\sigma_{c,\mathbf{y}}$ implies $\mathrm{Med}_c(\beta|\mathbf{y})$ to be less trustworthy and, hence, more evidence for $\beta = 0$; similarly and in fact more directly, the larger $O_{\mathbf{y}}$, the more evidence of shrinkage toward zero [indeed, when $O_{\mathbf{y}}$ exceeds $\Delta_c(\mathbf{y})$, we will shrink $\mathrm{Med}_c(\beta|\mathbf{y})$ all the way to zero!].

I provide this rather detailed derivation and discussion to demonstrate that achieving sparsity via the posterior median under a mixture prior is conceptually (and often computationally) more appealing than optimizing a penalized likelihood. With the latter approach, one appears to be using a black box, as determined by an optimization algorithm, much more so than with the former, where every step helps to gain probabilistic insights. Furthermore, one wonders if the $L^1$ aspect should go directly into the model, as with the pena-likelihood approach, or go into the loss function, as with the Bayesian approach (recall that median is the optimal estimator under the $L^1$ loss). But these concerns by no means undermine the value of the MLE approach; indeed, in Kong, Meng, Nicolae and Tan [7] and Kong, Meng and Nicolae [6], we identified a situation where it is logically impossible to perform any Bayesian inference, but MLE offers a very useful solution. Rather, they serve as an invitation to the authors to explore alternative methods for achieving sparsity, as well as to include the posterior median approach in their list of methods for comparison.

**Q.2.** *Should we really care about being only one-step, not multi-step?*

The authors clearly took pride in their one-step LLA estimator, a pride many authors would share. It achieves sparsity, requires least computation, and is "as efficient as the fully iterative method, provided that the initial estimators are reasonably good." Indeed, the authors' simulation results show—if I understand their notation correctly—that it is even capable of outperforming the estimator from the fully iterative algorithm it supposedly approximates. So why would anyone still care about doing more than one step?

I would, from both practical and theoretical points of view. Practically, the statement "provided that the initial estimators are reasonably good" worries me. Many methods work beautifully—at least one hopes—in the hands of experts, who have experiences, insight and the necessary skills to carry out whatever fine tuning is needed to achieve all the good properties as promised. But once it is "at large," a cute indoor pet might be too excited or too frightened to be as obedient as it is with its loving owner, especially when it sees a complete stranger. One advantage of the MM or EM-type algorithms, because of their ascent property, is their relative robustness to starting points, thereby reducing the burden of fine tuning for a novice. (Of course, this statement is true only when the objective function being maximized has a unique mode; but the multi-modality does not seem to be an issue addressed in the paper, and at any rate it is better handled by a full Bayesian method than by the MLE approach.) Other than for large simulation studies, the extra computation needed for going beyond one step or even for iterating until convergence is usually quite acceptable given current computing power. At any rate, we can view this extra computational cost as an insurance premium for guarding against accidental "unreasonable" starting points, assuming, of course, the fully iterated results are what we are really after.

This brings me to my second point. Theoretically, I am curious about how the ascent property helps to improve the performance of the $k$-step estimator, as $k$ increases. Or does it? If one considers higher (penalized) likelihood value as desirable, then of course the improvement is automatic. But does higher penalized likelihood value always transfer to better correct-fit rate? There is no evidence in the paper—or elsewhere as far as I am aware of—to suggest this. Quite to the contrary, the numerical evidence (e.g., Table 3) the authors provided seem to indicate that the opposite is more than a theoretical possibility, because there the one-step SCAD outperformed the fully iterated SCAD by a large margin.

This outperformance is rather unsettling for me. If higher penalized likelihood values do not necessarily transfer to better fitting proportion, which is one of the main criteria used, why then maximize the likelihood to start with? Could it be that the one-step LLA worked well not because it *approximates* the "ideal" full algorithm well, but rather it—perhaps by accident—*corrects* the defect in the full MLE? Actually, this is not the first time that such a phenomena has been noted. In Silverman, Jones, Wilson and Nychka [15], it was reported that stopping an EM iteration earlier can provide better image reconstructions than reaching the full MLE, due to the over-fitting tendency with the MLE approach, especially when there are many parameters involved. A remedy is to be more Bayesian: to rely on the entire posterior distribution rather than just its mode (even the use of the posterior median can help to reduce the over-fitting rate because of the robustness property of the median). One can even achieve this, at least partially, without paying the B-club membership; Silverman et al. [15] inserted a "smoothing step" into EM to regulate the over-fitting problem. I wonder if a similar idea could help to reduce the much higher over-fitting rates for the full SCAD, as compared to the one-step SCAD, in the authors' Table 3. Again, it would be useful to include the posterior median approach in the authors' comparisons, so we can have a better diagnosis as to whether the failure of the full SCAD is due to algorithmic inaccuracy or rather defects in the underlying principles or methods.

**3. Some ironies.**  As I am finishing this discussion, several ironies come to mind. First, usually a discussant is more critical of topics s/he has dived well into, and much more "generous" otherwise, in fear of making a fool of her/himself. I clearly have failed to follow this wisdom, as my "outsider's comments" touch upon more critical issues than my "insider's inspection" does. I hope that the authors will tolerate my adventurous style for wandering outside boundaries that I am comfortable with, and will take these outsider's comments only as seriously as the authors think they deserve to be.

Second, Section 2 of my discussion seems to complain loudly about the lack of Bayesian thinking, and practice, in the current paper or in much of the similar literature. But the authors' formulation of the "EM = MM" result is squarely Bayesian, as it invokes a prior and even a hyper-prior specification. In contrast, my original construction and derivation, as quoted in this discussion from Meng [10], is squarely likelihood-based. I want to mention this to emphasize the fact that—to echo the title—we are all multi-color cats and we all try to catch as many mice, or even rats, as we can.

Last, although I have tried to take my discussant's role very seriously (and perhaps too seriously since I just couldn't "let it go" even when it was three months overdue), I have not done the one thing that most discussants would do first, that is, to congratulate the authors for a stimulating paper! So here it comes: Congratulations and Meow!!!

**4. Acknowledgments.** I thank the Co-Editor, Jianqing Fan, for his kind invitation and his extraordinary patience in tolerating my tardiness in preparing this discussion, and Paul Baines, Xuming He and Thomas Lee for helpful comments.

Department of Statistics
Harvard University
Cambridge, Massachusetts 02138
USA
E-mail: meng@stat.harvard.edu